\documentclass[a4paper,11pt]{article}

\usepackage[english]{babel}
\usepackage{cite}

\usepackage{amsmath}
\usepackage{amssymb}
\usepackage{amsfonts}
\usepackage{array}
\usepackage{graphicx}

\usepackage{epsfig}
\usepackage{float}

\usepackage{latexsym}
\usepackage{amsthm}
\usepackage{enumerate}
\usepackage{geometry}
\renewcommand{\paragraph}[1]{\bigskip \textbf{#1}.}

\newcommand{\norm}[1]{\left\lVert#1\right\rVert}

\newcommand*{\partialfrac}[3]{\ensuremath{\frac{\partial^{#2} #1}{\partial #3^{#2}}}}

\newcommand{\de}{\partial}

\newcommand{\mrm}[1]{\mathrm{#1}}

\newcommand{\uvec}[1]{\boldsymbol{\textbf{#1}}}

\newcommand{\spaceV}{\mathcal{V}}
\newcommand{\spaceW}{\mathcal{W}}

\newcommand{\Tau}{\mathcal{T}}

\newcommand{\B}{\mathcal{B}}

\numberwithin{equation}{section}
\allowdisplaybreaks

\newtheorem{theorem}{Theorem}[section]

\newtheorem{remark}[theorem]{Remark}

\begin{document}
	
	\title{A diffuse interface box method for elliptic problems}

	\author{G. Negrini$^a$, N. Parolini$^a$ and M. Verani$^a$}

	\maketitle
	
	\begin{center}
		{\small
			$^a$ MOX, Dipartimento di Matematica, Politecnico di Milano, Piazza Leonardo da Vinci 32, I-20133 Milano, Italy
		}
	\end{center}

\begin{abstract}
	We introduce a {\em diffuse interface box method }(DIBM)  for the numerical approximation on complex geometries of elliptic problems with Dirichlet boundary conditions. We derive a priori $H^1$ and $L^2$ error estimates highlighting the r\^{o}le of the mesh discretization parameter and of the diffuse interface width. Finally, we present a numerical result assessing the theoretical findings. \\
	
\end{abstract}

\maketitle

\section{Introduction} \label{sec:intro}

The finite volume method (FVM) is a popular numerical
strategy for solving partial differential equations modelling real life problems. One crucial and attractive
property of FVM is that, by construction, many physical conservation laws possessed
in a given application are naturally preserved. Besides,
similar to the finite element method, the FVM can be used to deal with domains with complex geometries. In this respect, one crucial issue is the construction of the computational grid. To face this problem, one can basically resort to two different types of approaches. In the first approach, a mesh is constructed on a sufficiently accurate approximation of the exact physical domain (see, e.g.,  isoparametric finite elements \cite{Ciarlet:book}, isogeometric analysis \cite{IGA:book}, or Arbitrary Lagrangian-Eulerian formulation \cite{DGH,HAC,HLZ}), while in the second approach (see, e.g., Immersed Boundary methods \cite{Peskin:acta}, the Penalty Methods \cite{Babuska:penalty}, the Fictitious Domain/Embedding Domain Methods \cite{Borgers-Widlund:1990,BG03,BBSV:2019}, the cut element method \cite{BH10,BH12} and the Diffuse Interface Method \cite{Li-Lowengrub:2009}) one embeds the physical domain into a simpler computational mesh whose elements can intersect the boundary of the given domain. Clearly, the mesh generation process is extremely simplified in the second approach, while the imposition of boundary conditions requires extra work.  Among the methods sharing the second approach, in this paper we focus on the diffuse interface approach developed in \cite{pen:schlottbom}. In parallel, we consider, for its simplicity, the piecewise linear FVM, or {\em box method}, that has been the object of an intense study in the literature (see, e.g., the pioneering works \cite{box:bankrose, box:hackbusch} and the more recent \cite{box:ewinglin,box:xuzou}).

The goal of this paper is to propose and analyse a diffuse interface variant of the box method, in the sequel named DIBM (diffuse interface box method), obtaining a priori $H^1$ and $L^2$ error estimates depending both on the discretization parameter $h$ (dictating the accuracy of the approximation of the PDE) and the width $\epsilon$ of the diffuse interface (dictating the accuracy of the domain approximation). Up to our knowledge, this is new in the literature. Besides, the study of DIBM for elliptic problems, despite its simplicity, opens the door to the study of more complicated differential problems and to the analysis of diffuse interface variants of more sophisticated finite volume schemes.

The outline of the paper is as follows. In section \ref{sec:prelim} we briefly recall the box method, while in section \ref{sec:diffintbox}, we present the diffuse interface box method (DIBM) along with a priori error estimates.  Finally in section \ref{sec:experiments} we will provide a numerical test to validate the theoretical results. The numerical results have been obtained using the open-source library OpenFOAM{\sffamily\textregistered}.

\section{The box method} \label{sec:prelim}
In this section, we recall (see \cite{box:bankrose, box:hackbusch, box:xuzou}) the box method for the solution of an elliptic problem. Let $D \subset \mathbb{R}^2$ be a polygonal bounded domain (in the following section this hypothesis will be relaxed). We consider the following problem:
\begin{equation}
\begin{cases}
-\Delta u = f, \quad & \mrm{in}~ D \\
u=g, \quad & \mrm{on}~ \Gamma = \de D,
\end{cases}
\label{eq:strong:problem}
\end{equation}
where $f \in L^2(\Omega)$ and $g \in H^{1/2}(\Gamma)$. 

Let $\Tau_h = \{t_i\}$ be a conforming and shape regular triangulation of $D$. We denote by $h_t$ the diameter of $t \in \Tau_h$ and we introduce the set $\textsf{V}_h=\{\textsf{v}_i\}$ of vertices of $\Tau_h$ with $\textsf{V}_h = \textsf{V}_h^\de \cup \textsf{V}_h^o$, the set $\textsf{V}_h^o$ containing the interior vertices of $\Tau_h$. We denote by $w_\textsf{v}$ the set of triangles sharing the vertex $\textsf{v}$. On $\Tau_h$ we define the space of linear finite elements
$$
\spaceV_{h,g_h} = \left\lbrace v_h \in C^0(\bar{D}): v_h\vert_t \in \mathbb{P} ^1(t) ~\forall t \in \Tau_h ~\mrm{and}~ v_h = g_h ~\mrm{on}~ \de D\right\rbrace,
$$
where $g_h$ is a suitable piecewise linear approximation of $g$ on $\de D$.

Let $\B_h=\{b_\textsf{v}\}_{\textsf{v} \in {\textsf{V}}_h^o}$ be the ``box mesh'' (or dual mesh) associated to $\Tau_h$. Each box $b_\textsf{v}$ is a polygon with a boundary consisting of two straight lines in each related triangle $t \in w_\textsf{v}$. These lines are defined by the mid-points of the edges and the barycentres of the triangles in $w_\textsf{v}$.

On $\B_h$ we introduce the space of piecewise constant functions,
$$
\spaceW_h = \left\lbrace w_h \in L^2(D): w_h \in \mathbb{P}^0(b_\textsf{v}) ~\forall b_\textsf{v} \in \B_h \right\rbrace.
$$

The box method for the approximation of \eqref{eq:strong:problem} reads as follows: find $u_{B,h} \in \spaceV_{h,g_h}$ such that 
\begin{equation}
a_{\Tau_h}(u_{B,h},w_h)=(f,w_h)_D \quad \forall w_h \in \spaceW_h,
\label{eq:weak:problem}
\end{equation}
where
\begin{equation}
a_{\Tau_h}(u_{B,h},w_h) = - \sum_{\textsf{v} \in \textsf{V}_h^o} \int_{\de b_{\textsf{v}}}
\partialfrac{v_h}{}{\uvec{n}_b} w_h \mrm{d} s,
\label{eq:bilinearequivalence}
\end{equation}
being $\uvec{n}_b$ the outer normal to $b_\textsf{v}$ and $(\cdot,\cdot)_D$ is the usual $L^2$ scalar product on $D$.

Note that there holds (see \cite{box:bankrose, box:hackbusch} for the two dimensional case and \cite{box:xuzou} for the extension to any dimension)
\begin{equation}
\int_{\de b_{\textsf{v}}} \partialfrac{\phi_{\textsf{v}'}}{}{\uvec{n}_b} w_h \mrm{d} s
= \int_D \nabla \phi_\textsf{v} \cdot \nabla \phi_{\textsf{v}'} \mrm{d} x, \quad \forall \textsf{v} \in \textsf{V}_h^o, \forall \textsf{v}' \in \textsf{V}_h,
\label{eq:bilinearequivalenceboundary}
\end{equation}
where $\phi_{\textsf{v}}$ is the usual hat basis function with support equal to $w_\textsf{v}$.

The relation \eqref{eq:bilinearequivalenceboundary} is crucial to show the following perturbation results (see \cite{box:hackbusch, box:xuzou}):
\begin{equation}
\begin{aligned}
\norm{\nabla(u_{B,h}- u_{G,h})}_{L^2(D)} &\leq C h \norm{f}_{L^2(D)}, \\
\norm{u_{B,h}- u_{G,h}}_{L^2(D)} &\leq C h^2 \norm{f}_{L^2(D)}, \\
\end{aligned}
\label{eq:perturbationboxfem}
\end{equation}
where $h= \max_{t \in \Tau_h} h_t$ and $u_{G,h} \in \spaceV_{h,g_h}$ is the linear finite element approximation to the solution of problem \eqref{eq:strong:problem}.

\section{The box method with diffuse interface (DIBM)} \label{sec:diffintbox}
\begin{figure}[t]
	\centering
	\includegraphics[width=0.4\linewidth]{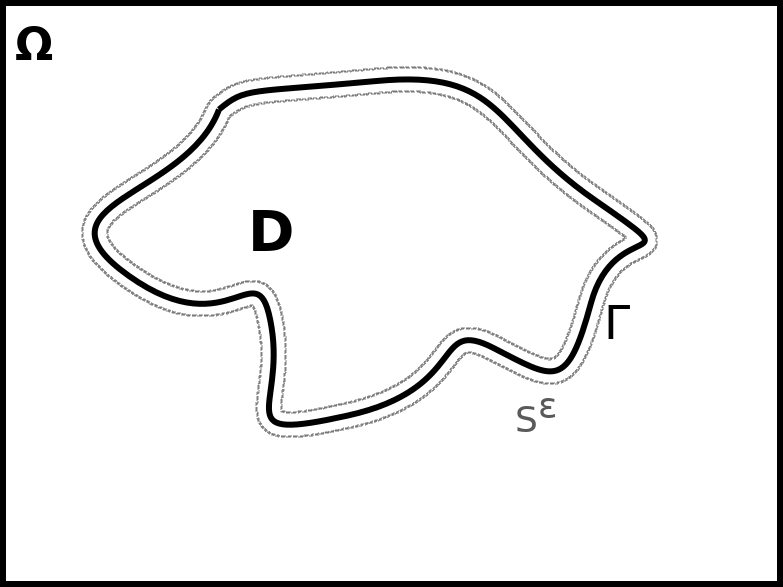}
	\caption{Diffuse interface representation: $D$ is a surrogate domain of $\Omega$; $\Gamma$ is the Dirichlet boundary and $S^\epsilon$ is its tubular neighbour.}
	\label{fig:domain:diffint}
\end{figure}

The aim of this section is to introduce a variant of the box method for the approximate solution of problem \eqref{eq:strong:problem} in case of a general (non-polygonal) domain $D \subset \mathbb{R}^2$, where in the spirit of \cite{pen:schlottbom} the Dirichlet boundary condition is treated with a diffuse interface approach. To this aim we introduce an hold-all domain $\Omega$ such that $D \subset \Omega$. In the sequel we will work under the hypothesis $\Gamma = \de D \in C^{1,1}$. With a slight abuse of notation we denote by $\Tau_h$ a shape regular triangulation of $\Omega$. It is worth noting that $\Tau_h$ is not conforming with $D$. Following \cite{pen:schlottbom} we first select a tubular neighbourhood $S^\epsilon$ of $\Gamma$, where $\epsilon$ denotes the width of $S^\epsilon$ (see Figure \ref{fig:domain:diffint}). Then we introduce the set $S^\epsilon_h$ which contain s all the triangles of $\Tau_h$ having non-empty intersection with $S^\epsilon$. Note that the width of the discrete tubular neighbourhood $S^\epsilon_h$ is $\delta+ \epsilon$ where $\delta$ is the maximum diameter of triangles crossed by $\de S^\epsilon$.

To proceed, we assume that there exists an extension $\tilde{g} \in H^2(\Omega)$ of the boundary data g.

We set $D^\epsilon_h = D \backslash S^\epsilon_h$ and introduce the function $u^{\epsilon,h} \in H^1(D_h^\epsilon)$ such that $u^{\epsilon,h}=g$ on $\de D^\epsilon_h$, which solves the following continuos problem:
\begin{equation}
\int_{D^\epsilon_h} \nabla u^{\epsilon,h} \cdot \nabla v = \int_{D^\epsilon_h} f v \quad \forall v \in H^1_0(D^\epsilon_h).
\label{eq:weak:diffintcontinuous}
\end{equation}

The solution $u^{\epsilon,h}$ is then extended to $S^\epsilon_h$ by setting $u^{\epsilon,h}=\tilde{g}$ in $S^\epsilon_h$.

The following results have been proved in \cite[Thm 1.2]{pen:schlottbom}:
\begin{equation}
\frac{1}{\epsilon + \delta} \norm{u - u^{\epsilon,h}}_{L^2(D)} +
\frac{1}{\sqrt{\epsilon + \delta}} \norm{\nabla u - \nabla u^{\epsilon,h}}_{L^2(D)} 
\leq C \left( \norm{f}_{L^2(D)} + \norm{g}_{H^2(D)} \right).
\label{eq:perturbation:diffint}
\end{equation}

Let $\spaceV_{h,\tilde{g}_h}^\epsilon = \left\lbrace
v_h \vert_{D^\epsilon_h}: v_h \in \mathbb{P}^1(t) \forall t \in \Tau_h ~\mrm{and}~v_h = \tilde{g}_h ~\mrm{on}~ \de D^\epsilon_h \right\rbrace$, with $\tilde{g}_h$ the Lagrangian piecewise linear interpolant of $\tilde{g}$.

It has been proved (cf. \cite[Thms 5.1 and 5.3]{pen:schlottbom}) that the linear finite element approximation $u^\epsilon_{G,h} \in \spaceV_{h,\tilde{g}_h}^\epsilon$ of $u^{\epsilon,h}$ satisfies the following estimates:
\begin{equation}
\begin{aligned}
\norm{\nabla (u^{\epsilon,h} - u^\epsilon_{G,h})}_{L^2(D)} &\leq
C (\sqrt{\delta} + \kappa^{\frac{2}{3}} + h) \left( \norm{f}_{L^2(D)} + \norm{\tilde{g}}_{H^2(D)} \right), \\
\norm{u^{\epsilon,h} - u^\epsilon_{G,h}}_{L^2(D)} &\leq
C (\delta + \kappa^{\frac{4}{3}} + h^2) \left( \norm{f}_{L^2(D)} + \norm{\tilde{g}}_{H^2(D)} \right), \\
\end{aligned}
\label{eq:estimate:diffint1}
\end{equation}
where $\kappa$ is the maximum diameter of the triangles intersection $\de S^{\epsilon + h}$ and $u^\epsilon_{G,h}$ has been extended to $D^\epsilon_h$ by setting $u^\epsilon_{G,h} = \tilde{g}_h$ on $S^\epsilon_h$. 

Let us now introduce the box method with diffuse interface (DIBM). We denote by $u^\epsilon_{B,h} \in \spaceV_{h,\tilde{g}_h}^\epsilon$, the approximation obtained from applying the box method to \eqref{eq:weak:diffintcontinuous} (cf. \eqref{eq:weak:problem}). The solution $u^\epsilon_{B,h}$ is then extended to $D$ by setting $u^\epsilon_{B,h}=\tilde{g}_h $ in $S^\epsilon_h$. Then employing the triangle inequality in combination with \eqref{eq:perturbation:diffint}, \eqref{eq:estimate:diffint1} and \eqref{eq:perturbationboxfem} we get the following estimates for DIBM:
\begin{equation}
\begin{aligned}
\norm{\nabla(u - u^\epsilon_{B,h})}_{L^2(D)} &\lesssim \sqrt{\epsilon + \delta} + \sqrt{\delta} + k^{\frac{2}{3}} + h, \\
\norm{u - u^\epsilon_{B,h}}_{L^2(D)} &\lesssim \epsilon + \delta+ k^{\frac{4}{3}} + h^2.
\end{aligned}
\label{eq:estimate:final}
\end{equation}

\begin{figure}[t]
	\centering
	\includegraphics[width=0.4\linewidth]{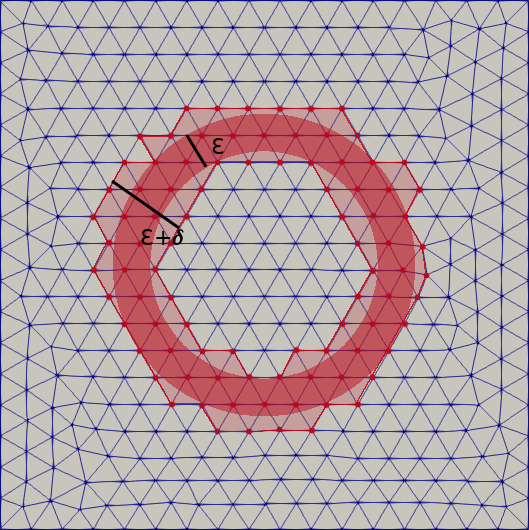}
	\includegraphics[width=0.4\linewidth]{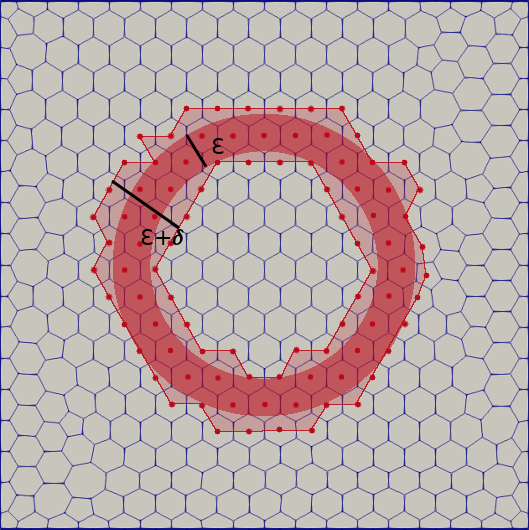}
	\caption{Discrete diffuse interface representation on triangulation (left) and on box mesh (right). Constrained cells are marked with red dots while the continuous and discrete diffuse interfaces are coloured by darker an lighter red respectively. }
	\label{fig:diffint}
\end{figure}

\section{Numerical experiments} \label{sec:experiments}
In this section we numerically assess the theoretical estimates obtained in Section \ref{sec:diffintbox}. To this aim, we consider the test case originally introduced in \cite[Section 6]{pen:schlottbom} that is briefly recalled in the sequel. Let $\Omega=(-1,1)^2$ and let $\Gamma$ be the boundary of the circle $B_1(0)$ with centre $(0,0)$ and unitary radius. Thus, $\Gamma$ splits the domain $\Omega$ into two subregions: $D_1=B_1(0)$ and 
$D_2=\Omega\setminus \overline{D}_1$. 
Let $u$ be the solution of the following problem 
\begin{equation}
-\Delta u = f ~~\text{in~}\Omega,\qquad u=g ~~\text{on~}\Gamma,\qquad
u=0 ~~\text{on~}\partial\Omega,
\end{equation}
where $g(x,y) = (4-x^2)(4-y^2)$ on $\Gamma$ and extended to $\Omega$ as
$
\tilde{g}(x,y)=(4-x^2)(4-y^2)\cos(1-x^2-y^2).
$

Setting the solution  equal to:
\begin{equation}
u(x,y) = (4-x^2)(4-y^2)\left(\chi_{D_2} + \exp(1- x^2 - y^2)\chi_{\bar{D}_1}\right),
\label{eq:analyticsol}
\end{equation}
where $\chi_{D_i}$, $i=1,2$ are the characteristic functions of the two parts of $\Omega$,  the source term  $f$ is chosen as:
\begin{equation*}
f=
\begin{cases}
-\Delta u &\quad \mrm{in}~ \Omega\backslash\Gamma,\\
0 &\quad\mrm{on}~\Gamma.
\end{cases}
\end{equation*}

All the computations have been performed employing a Voronoi dual mesh of a Delaunay triangulation (i.e., the dual mesh is obtained by connecting  the barycentres of the triangles with straight lines).

To validate the estimates \eqref{eq:estimate:final} we consider in a separate way the influence of $h$ and $\epsilon$ on the error. More precisely, we first explore the convergence with respect to $h$ and then we study the convergence with respect to $\epsilon$. In both cases we consider  a uniform discretization of the domain $\Omega$ so to have $\kappa=\delta=h$.

\paragraph{Convergence w.r.t. $\boldsymbol{h}$} We set $\epsilon=2^{-20} \ll h$ while we let $h$ vary as
$$
h = 0.056, 0.028, 0.0139, 0.00694.$$

\begin{figure}[ht]
	\centering
	\includegraphics[width=0.63\linewidth]{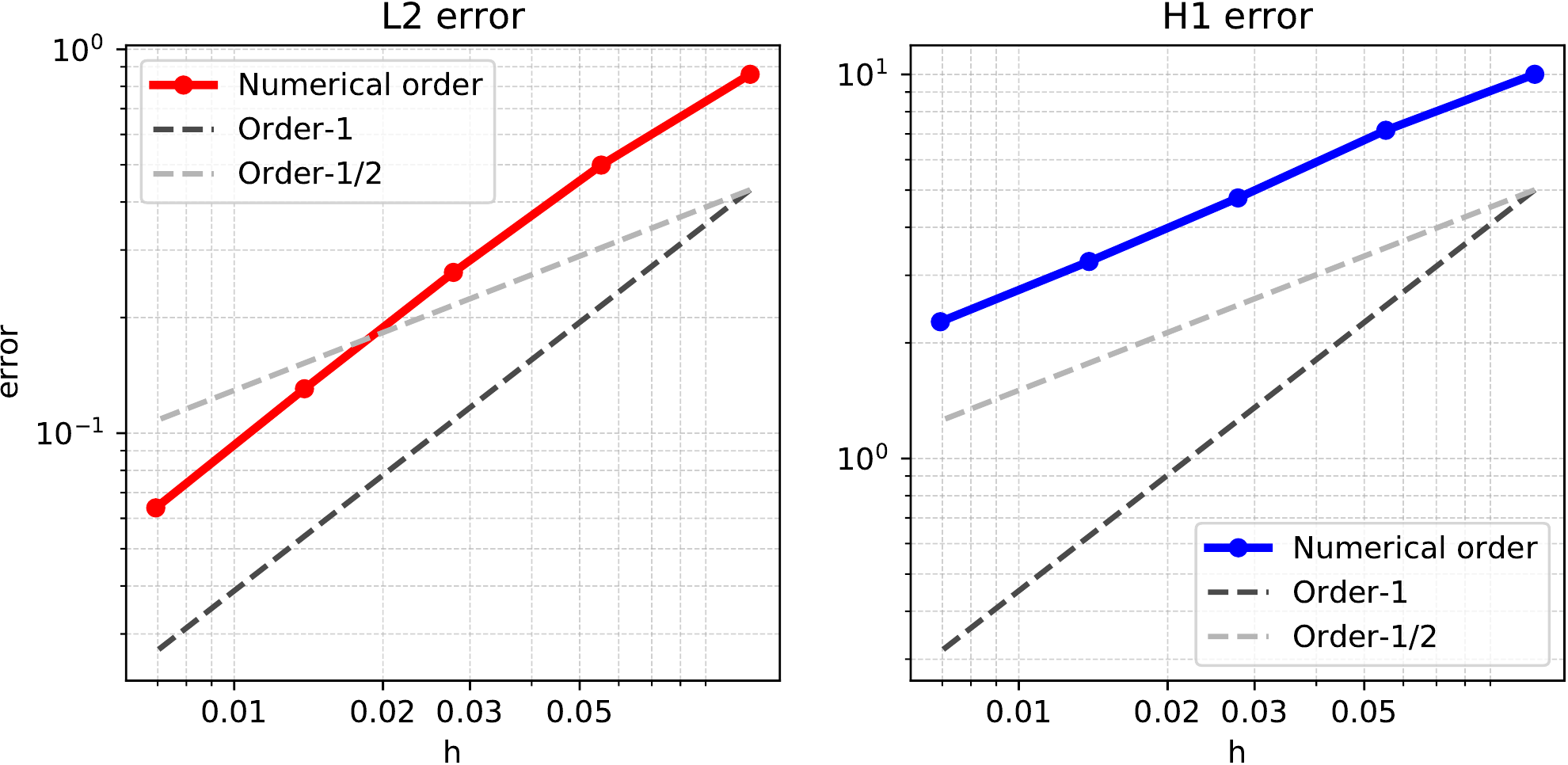}
	\caption{Error behaviour with respect to $h$ (fixed $\epsilon=2^{-20}$): (left) $L^2$-norm error, (right) $H^1$-norm error. Dashed lines are theoretical convergence orders.}
	\label{fig:diffint:convergenceH}
\end{figure}

From Figure \ref{fig:diffint:convergenceH} we observe that the $L^2$-norm of the error decreases with order $1$ while the error decreases with order 1/2 in the $H^1$-norm. These rates of convergence are in agreement with  \eqref{eq:estimate:final}.

\begin{figure}[ht]
	\centering
	\includegraphics[width=0.3\linewidth]{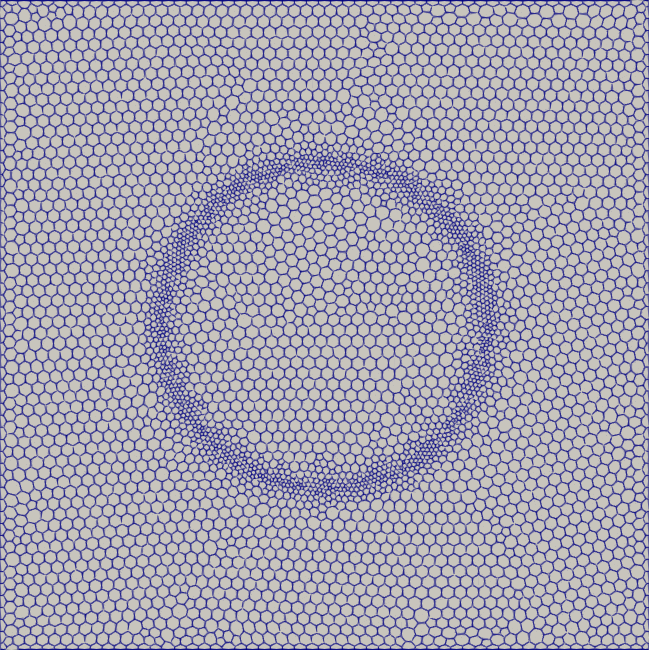}
	~~~~~
	\includegraphics[width=0.6\linewidth]{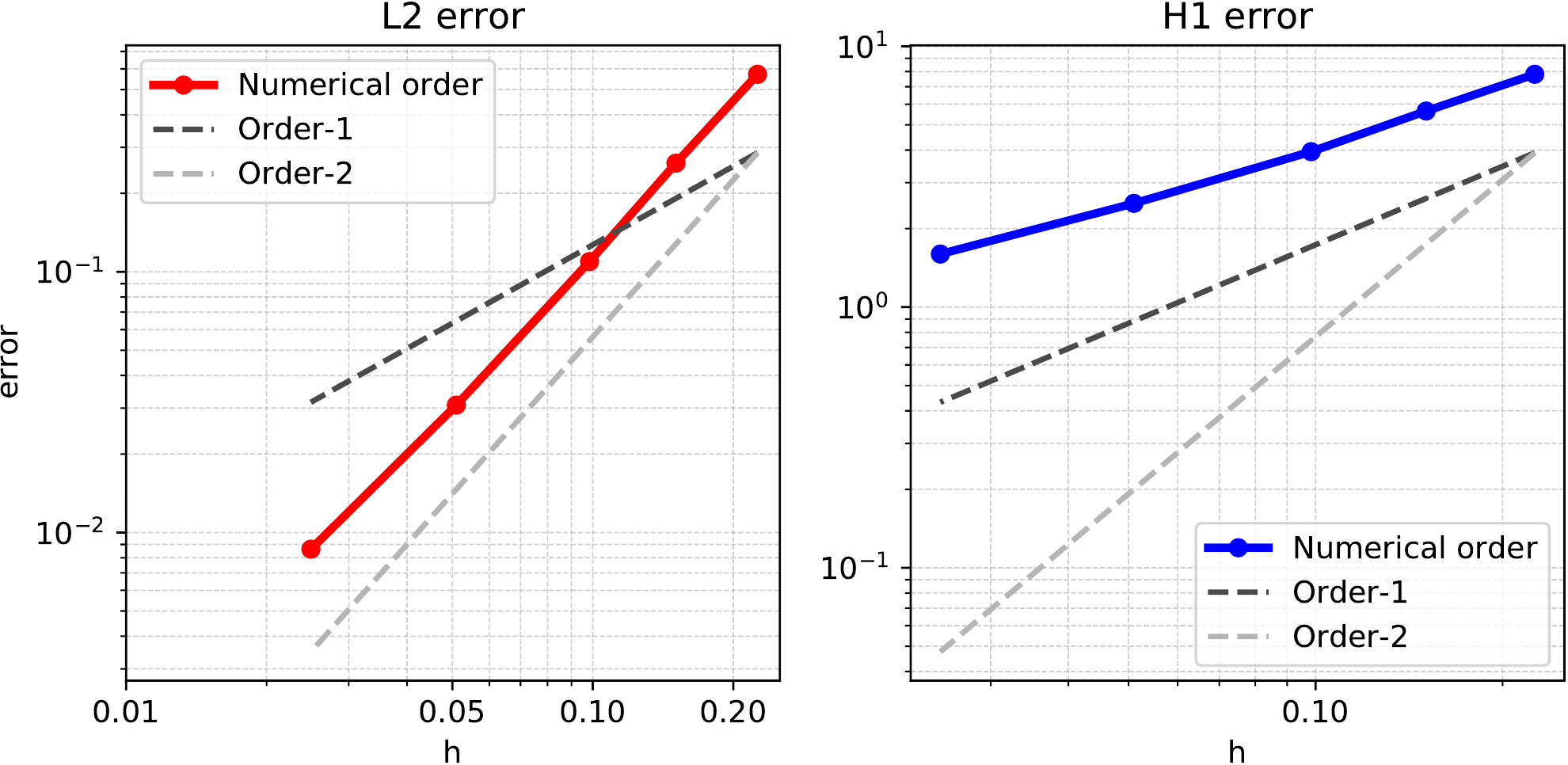}
	\caption{On the left: example of a dual mesh with local mesh refinement around surrogate boundary. On the right: error behaviour with respect to $h$ with local mesh refinement around the interface (fixed $\epsilon=2^{-20}$): (left) $L^2$-norm error, (right) $H^1$-norm error. Dashed lines are theoretical convergence orders.}
	\label{fig:diffint:convergenceHref}
\end{figure}

\begin{remark}
	If a local refinement of the diffuse interface region  is performed in such a way that $\delta \simeq \kappa \simeq h^2$ (Figure \ref{fig:diffint:convergenceHref}), then first and second order of convergence are recovered for $H^1$ and $L^2$ norms, respectively (cf. \cite[Section 6]{pen:schlottbom}).
\end{remark}

\begin{figure}[ht]
	\centering
	\includegraphics[width=0.63\linewidth]{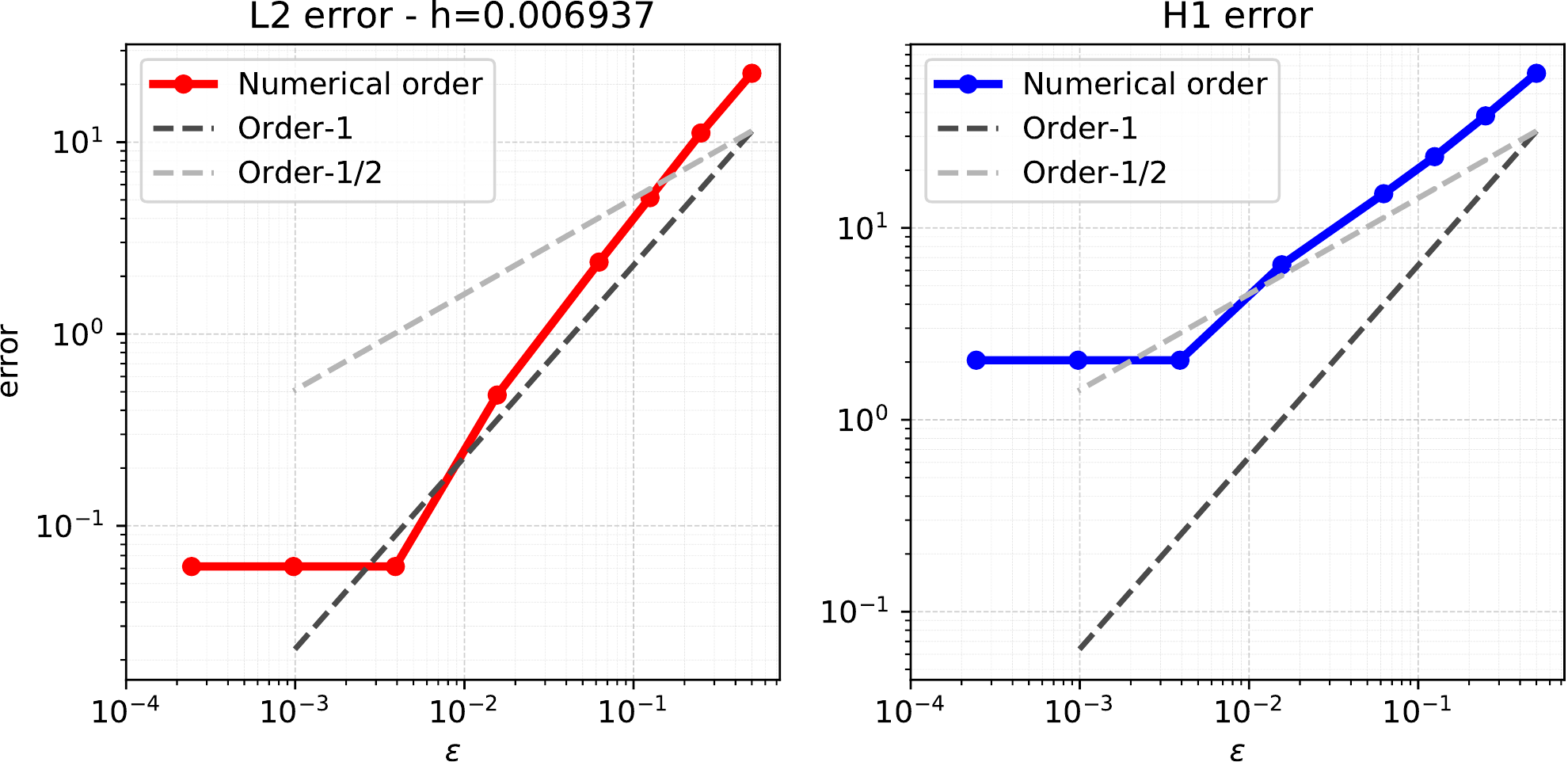}
	\caption{Error behaviour with respect to $\epsilon$ (fixed $h=0.00694$): (left) $L^2$-norm error, (right) $H^1$-norm error. Dotted lines are theoretical convergence orders.}
	\label{fig:diffint:convergenceE}
\end{figure}

\paragraph{Convergence w.r.t. $\boldsymbol{\epsilon}$} We employ a fine mesh ($h=0.00694$) and let the value of $\epsilon$ vary as:
$$
\epsilon = 2^{i},~i=-1,...,-20.
$$
The results are collected in Figure \ref{fig:diffint:convergenceE}. The theoretical rates of convergence  with respect to $\epsilon$ (cf. \eqref{eq:estimate:final}) are obtained both in the $L^2$-norm (order $1$) and in the $H^1$- norm (order $1/2$ ). It is worth noticing that when the value of $\epsilon$ becomes smaller than the chosen value of $h$, a plateau is observed as the (fixed) contribution from the discretization of the PDE (related to $h$) dominates over the contribution from  the introduction of the diffuse interface (related to $\epsilon$).

\section{Conclusions} \label{sec:conclu}
In this paper we introduced a diffuse interface variant of a finite volume method, namely of the the so-called {\em box method} and obtained $L^2$ and $H^1$ error estimates highlighting the contributions from the discretization parameter $h$ associated to the polygonal computational mesh and the width $\epsilon$ of the diffuse interface. Despite the  simplicity of the method, the present contribution seems to be novel in the literature. Moreover, the present work may represent the first step towards the study of the diffuse interface variant  of more sophisticated finite volume schemes (possibly for more complex differential problems).

This work opens fictitious boundary methods analysis to the box method and finite volume framework. Possible extensions of this research could be to being able to apply the plenty of penalization methods that are mostly thought for finite element implementations such as shifted boundary, Nitsche penalty, cut-fem or Brinkman penalization.

\section{Acknowledgements} \label{sec:aknowled}

The first author acknowledges the financial support of Fondazione Politecnico. The third author acknowledges the financial support of PRIN
research grant number 201744KLJL ``\emph{Virtual Element Methods:
	Analysis and Applications}'' funded by MIUR. The second and third authors acknowledge the financial support of INdAM-GNCS.


\bibliography{biblio}
\nocite{*}

\end{document}